\documentclass[letterpaper, 10 pt, conference]{ieeeconf}  

\IEEEoverridecommandlockouts                              
\usepackage{graphics}
\usepackage{graphicx} 
\usepackage{rotating}
\usepackage{dcolumn}
\usepackage{longtable}
\usepackage{multirow}
\usepackage{amsmath}
\usepackage{amssymb}
\usepackage{xspace}
\usepackage{textcase}
\usepackage[margin=1in]{geometry}
\usepackage{float}
\usepackage[latin1]{inputenc}
\usepackage{tikz}
\usetikzlibrary{trees}
\usetikzlibrary{arrows}

\usepackage{xcolor} 
\usepackage{ulem} 

\newtheorem{theorem}{Theorem}[section]
\newtheorem{corollary}{Corollary}[theorem]

\newtheorem{definition}{Definition}[section]

\title{\LARGE \bf
Koopman Representations of Dynamic Systems with Control
}

\author{Craig Bakker, W. Steven Rosenthal, and Kathleen E. Nowak
\thanks{This work was supported by project [] under the Deep Science Agile Investment at PNNL}
\thanks{C. Bakker, K. Nowak, and W.S. Rosenthal are with the Pacific Northwest National Laboratory,
        Richland, Washington
        {\tt\small craig.bakker@pnnl.gov}}%
} 

\begin{document}

\maketitle
\thispagestyle{empty}
\pagestyle{empty}

\begin{abstract}

The design and analysis of optimal control policies for dynamical systems can be complicated by nonlinear dependence in the state variables. Koopman operators have been used to simplify the analysis of dynamical systems by mapping the flow of the system onto a space of observables where the dynamics are linear (and possibly infinte). This paper focuses on the development of consistent Koopman representations for controlled dynamical system. We introduce the concept of dynamical consistency for Koopman representations and analyze several existing and proposed representations deriving necessary constraints on the dynamical system, observables, and Koopman operators. Our main result is a hybrid formulation which independently and jointly observes the state and control inputs. This formulation admits a relatively large space of dynamical systems compared to earlier formulations while keeping the Koopman operator independent of the state and control inputs. More generally, this work provides an analysis framework to evaluate and rank proposed simplifications to the general Koopman representation for controlled dynamical systems.

\end{abstract}

\section{INTRODUCTION}


The Koopman operator provides a way to transform a (potentially) nonlinear finite-dimensional dynamical system into an infinite-dimensional linear system.  It does this by lifting the nonlinear state dynamics into a functional space of observables, where the dynamics are linear \cite{budisic12jsr}.  Koopman representations of nonlinear systems are typically calculated using data-driven methods -- in particular, the use of time series data to calculate finite truncations of the Koopman operator and its associated observables.  For an example of an analytical Koopman representation, though, see Page and Kerswell \cite{page18jsr}.

Dynamic Mode Decomposition (DMD) \cite{tu14jsr} is a common technique for doing this.  DMD works by defining Koopman eigenfunctions as linear combinations of state variable measurements.  Extended DMD (EDMD) builds on this by working with nonlinear functions of the state measurements \cite{williams15jsr}. EDMD has the ability to represent Koopman observables (or, equivalently, Koopman eigenfunctions) that are nonlinear functions of the state.  This can result in greater accuracy, but it also then creates the challenge of choosing a good dictionary of functions from which to work \cite{williams15jsr,kutz16jsr}.  Common choices include sets of polynomials \cite{williams15jsr,proctor18jsr} and radial basis functions \cite{williams16cp,korda18jsr}.  The expansions of the basis as the size of the state space increases can quickly become unmanageable, though, as an outworking of the curse of dimensionality.  

Williams et al. proposed using a kernel method and demonstrated some success with it \cite{williams14jsr} in tackling the basis expansion problem.  Other researchers have chosen to use neural networks to learn the observable functions from the input data rather than selecting a function dictionary \textit{a priori} \cite{yeung17jsr,li17jsr}.  As a nonlinear regression, this approach requires more computational effort to find the Koopman approximation, but it offers greater flexibility while requiring less initial insight into the problem.  Moreover, Yeung et al. were able to get significantly higher accuracy than EDMD using polynomials \cite{yeung17jsr}.

The ability to represent nonlinear dynamical systems as linear systems in a non-local way can be valuable from a control perspective, and indeed several researchers have started to look at this.  Brunton et al. start by considering a small system with an analytical finite Koopman representation \cite{brunton16jsr}.  They then use a Linear Quadratic Regulator (LQR) formulation to get a feedback control law.  This approach performs significantly better than LQR applied to local state linearizations.  Small changes in the underlying dynamical system, however, produce a control system that is no longer Koopman invariant, and this creates problems for applying LQR.

Proctor et al. \cite{proctor18jsr,proctor16jsr} show how to use DMD when there are control inputs applied to the system.  This is essentially just system identification within the Koopman framework -- the paper does not consider the formulation of a control policy, optimal or otherwise, from the data.  The formulation used assumes an affine control (i.e., $B u$) rather than lifting the control inputs to the space of observables, but given that DMD only produces linear combinations of input data, using $B u$ and solving for $B$ is essentially equivalent to this.  Williams et al. \cite{williams16cp} provides the same contribution but for EDMD rather than DMD, and their approach does modify the control inputs in a nonlinear way.  

Unlike most Koopman control papers, which work in discrete time, Kaiser et al. \cite{kaiser17jsr} use a continuous-time Koopman representation, and they then apply LQR to that.  Their formulation produces an affine control form with non-constant coefficients for the control term.  They are still able to solve for a feedback control law using the appropriate form of the Riccati equation.  Korda and Mezi\'{c} \cite{korda18jsr} work with a discrete-time controller (applied to a continuous-time system via discretization), and they combine the Koopman representation with Model Predictive Control (MPC) and a receding time horizon approach.  Using MPC denies them the possibility of a closed-form feedback control law, but it does allow them to include additional constraints that would not be present in the LQR formulation, and the optimization problem being solved is still convex.  The combination offline calculations to calculate the truncated Koopman representation and to reduce the problem dimension by using the problem's `dense' form result in a very fast online solver that performs noticeably better than a local linear model.  Arbabi et al. employ a similar approach \cite{arbabi18jsr}.

In this paper, we propose a method for analyzing Koopman representations of controlled dynamical systems. Section \ref{sec:KoopRep} provides some background on Koopman operators, while Section \ref{sec:consist} introduces necessary constraints for a given representation to be able to capture the flow of discrete- and continuous-time dynamical system with and without control inputs. We analyze several proposed Koopman representations under this framework. Section \ref{discussion} discusses the implications of our results in terms of the generality of controlled dynamical systems that can represented by these Koopman formulations.  We provide some brief concluding remarks in our conclusion.

\section{Koopman Operators for Dynamical Systems\label{sec:KoopRep}}

For a system whose dynamics are determined by both a state and control vector, a set of functions called \textit{observables} encode the state and control variables at a given time into a vector (of possibly infinte length). The Koopman operator propagates this information forward in time, still in terms of the set of observables. If the Koopman operator and observables can be derived or computed for a given dynamical system, then the linear form of the Koopman operator may make it easier (or possible) to analyze the dynamical system, as opposed to the original nonlinear formulation. However, the observables and the Koopman operator can be nonlinear functions of the state and control variables.

We provide the Koopman operator first for autonomous dynamics and then for controlled systems. Consider a dynamical system $\dot{x} = f\left(x\right), \ x \in M \subseteq \Re^n$ and a (possibly infinite) set of observables $\phi_j:M \rightarrow \Re$.  The flow mapping $F^{\tau}\left(x\left(t\right)\right) = x\left(t + \tau\right), \ F:M\rightarrow M$ maps points along the flow defined by $f$.  The Koopman operator $\mathcal{K}^{\tau}$ is defined as
\begin{gather}
\mathcal{K}^{\tau} \phi_j = \phi_j \circ F_{\tau} \\
\mathcal{K}^{\tau} : \mathcal{F} \rightarrow \mathcal{F} \\
\phi_j \circ F^{\tau} : M \rightarrow \Re
\end{gather}
\noindent where $\phi_j$ are functions existing in a functional space $\mathcal{F}$.  The set of $\phi_j$ functions, $S_{\phi}$, is invariant under the Koopman operator if the image of $S_{\phi}$, under the Koopman operator, is contained within the span of $S_{\phi}$ -- i.e., if $\mathcal{K}^{\tau} \phi_j$ is a linear combination of the elements of $S_{\phi}$ for all $\phi_j \in S_{\phi}$. Hereafter, the superscript $\tau$ may be dropped for notational simplicity.

In the case of a continuously differentiable set of observables $\phi_j$, then it makes sense to consider the instantaneous ($\tau\to 0$) flow induced by the Koopman operator in the space of observables,
\begin{gather}
\partial_t \phi = \mathcal{K}^{\tau=0^-} \phi = \mathcal{L}\phi = f\left(x\right)\cdot \nabla x
\end{gather}
Here, $\mathcal{L}$ is the infinitesimal generator for the Koopman semi-group, and provides the theoretical basis for the application of the chain rule to each observable $\phi[x(t)]$.  When the (possibly infinite) basis $\psi = \{\phi_1, ..., \phi_{N_\phi}, ...\}$ for $S_{\phi}$ needs to be considered, $\psi$ may be used in place of $\phi$. In practice, $\psi$ may be truncated to make numerical computation practical, in which case $S_{\phi}$ may be only approximatly invariant. In some cases it is advantagous to require the Koopman operator to reproduce the state-space dynamics in addition to the flow on $S_\phi$, such as in the Korda-Mezi\'{c} and Arbabi formulations \cite{korda18jsr, arbabi18jsr}.  In this case, $\psi$ contains $\mbox{Id}[x] = x$, and the set of observables is said to be \textit{state-inclusive}.  Intuitively, it makes sense to assume that observing the state directly is beneficial for a control policy. In a later section, this assumption will impact dynamic consistency and motivate Koopman formulations that jointly observe the state and control inputs.

\section{Consistent Koopman Representations for Dynamical Systems with Control \label{sec:consist}}

In the most general context, the Koopman representation of a dynamical system that depends on both state and control variables is
\begin{gather}
\psi\! \left( x_{k+1}, u_{k+1} \right) = \mathcal{K}\! \left( x_k, u_k \right) \psi\!\left( x_k, u_k \right) \\
\partial_t \psi\! \left( x, u \right) = \mathcal{L}\! \left( x, u \right) \psi\!\left( x, u \right)
\end{gather}

\noindent for discrete- and continuous-time, respectively.  The Koopman operator defines a new dynamical flow in a space of observables, e.g. $S_{\phi}$. The ability of a proposed Koopman representation to represent the dynamics of the underlying system will depend on the form of the representation and the nature of the underlying system. In this section, we derive necessary conditions for a Koopman representation to be consistent with the underlying system and evaluate several formulations used in the literature with respect to these conditions.

Williams et al. \cite{williams16cp} use
\begin{gather}
\psi\left(x_{k+1}\right) = K\left(u_k\right) \psi \left(x_k\right) \label{eqn:will_form} \\
K\left(u\right) = \sum_i \psi^u_i \left(u\right) K_i
\end{gather}
\noindent in which the observables depend only on the state, while the flow induced by the Koopman operator can be a function of the control input. Both Korda and Mezi\'{c} \cite{korda18jsr} and Arbabi et al. \cite{arbabi18jsr} do not lift the controls to the space of observables, and propose
\begin{equation} \label{eqn:kma_form}
\psi\left(x_{k+1}\right) = K \psi \left(x_k\right) + B u_k
\end{equation}
\noindent an affine control problem with constant control coefficients.  Kaiser et al. \cite{kaiser17jsr} use a joint observable of the form

\begin{gather}
\partial_t \psi \left(x,u\right) = \Lambda \psi \left(x,u\right) + \frac{\partial \psi}{\partial u} \dot{u} \label{eqn:kaiser_form}
\end{gather}

\noindent where $\Lambda = \text{diag} \left(\lambda_1,\lambda_2,\ldots\right)$ is a diagonal matrix of Koopman eigenvalues and the observerables are Koopman eigenfunctions.

For the joint state $(x,u)$, the general Koopman formulation will admit any nonlinear dynamical system. When simplifications like the above are introduced -- say to be tractable or admit an analytical solution -- then the simpler formulation may not be consistent with every dynamical system.  This inconsistency would then result in an unresolvable bias or model error.

\subsection{Dynamical consistency \label{subsec:cond}}

To illuminate the differences between these and and later Koopman representations, the following definitions propose a consistency requirement which determines whether a Koopman representation can represent a given discrete or continuous dynamical system.
\vspace{.5pc}
\begin{definition}
A Koopman representation $(\mathcal{L}, \psi)$ of a continuous-time dynamical system $\dot{x} = f\left(x\right)$ is \textit{consistent} with that system if the observables $\psi\!\left(x\right)$ satisfy
\begin{gather}
\frac{\partial \psi}{\partial t} = \frac{\partial \psi}{\partial x} f\left(x\right)
\end{gather}
and a Koopman representation of continuous-time system with control $\dot{x} = f\left(x,u\right)$ is \textit{consistent} with that system if $\psi\! \left(x\right)$ satisfy
\begin{gather}
\frac{\partial \psi}{\partial t} = \frac{\partial \psi}{\partial x} f\left(x,u\right)
\end{gather}

If the observable is a joint observable $\psi\left(x,u\right)$, then 

\begin{gather}
\frac{\partial \psi}{\partial t} = \frac{\partial \psi}{\partial x} f\left(x,u\right) + \frac{\partial \psi}{\partial u} \dot{u} \label{eqn:joint_consist}
\end{gather}

\end{definition}
\vspace{.5pc}
\begin{definition}
A Koopman representation $(\mathcal{K},\psi)$ of a discrete-time dynamical system $x_{k+1} = f\left(x_k\right)$ is \textit{consistent} with that system if the observables $\psi\left(x\right)$ satisfy
\begin{gather}
\frac{\partial \psi_{k+1}}{\partial x_k} = \frac{\partial \psi_{k+1}}{\partial x_{k+1}} \frac{\partial f_k}{\partial x_k}
\end{gather}
where $\left(\cdot\right)_k = \left. \left(\cdot\right) \right|_{x_k}$, and a Koopman representation of a discrete-time system with control $x_{k+1} = f\left(x_k,u_k\right)$ is \textit{consistent} with that system if $\psi \left(x\right)$ satisfy
\begin{gather}
\frac{\partial \psi_{k+1}}{\partial x_k} = \frac{\partial \psi_{k+1}}{\partial x_{k+1}} \frac{\partial f_k}{\partial x_k} \\
\frac{\partial \psi_{k+1}}{\partial u_k} = \frac{\partial \psi_{k+1}}{\partial x_{k+1}} \frac{\partial f_k}{\partial u_k}
\end{gather}

A consistent joint observable $\psi\left(x,u\right)$ would satisfy

\begin{gather}
\frac{\partial \psi_{k+1}}{\partial x_k} = \frac{\partial \psi_{k+1}}{\partial x_{k+1}} \frac{\partial f_k}{\partial x_k} \\
\frac{\partial \psi_{k+1}}{\partial u_k} = \frac{\partial \psi_{k+1}}{\partial x_{k+1}} \frac{\partial f_k}{\partial u_k} + \frac{\partial \psi_{k+1}}{\partial u_{k+1}} \frac{\partial u_{k+1}}{\partial u_k}
\end{gather}

\end{definition}
\vspace{.5pc}
These definitions describe necessary (though possibly not sufficient) conditions for a dynamical system to have a particular Koopman representation, as shown in the following theorem.
\vspace{.5pc}
\begin{theorem} \label{thm:dynconsist}
If $(\mathcal{K}, \psi)$ or $(\mathcal{L}, \psi)$ is a Koopman representation of a discrete-time or continuous-time dynamical system, respectively, then it is consistent with that dynamical system.
\end{theorem}
\vspace{.5pc}
\begin{proof}
Suppose $(\mathcal{L}, \psi)$ is a Koopman representation for a continuous-time dynamical system. The result follows directly from the chain rule applied to each observable $\phi\left[ x(t) \right] \in \psi$. If $\dot{x} = f\left(x\right)$, then we have
\begin{gather}
\frac{\partial }{\partial t} \left(\psi\left(x\right)\right) = \frac{\partial \psi}{\partial x} \frac{\partial x}{\partial t} = \frac{\partial \psi}{\partial x} f(x)
\end{gather}
The results for a controlled dynamical system $\dot{x} = f\left( x, u \right)$ follow analogously.  Similarly, suppose $(\mathcal{K}, \psi)$ is a Koopman representation for the discrete-time dynamical system, $x_{k+1} = f\left(x_k\right)$. Then,
\begin{gather}
\psi \left(x_{k+1}\right) = \psi \left(f\left(x_k\right)\right) \\
\frac{\partial }{\partial x_k} \left( \psi \left(x_{k+1}\right) \right) = \frac{\partial \psi_{k+1}}{\partial x_{k+1}} \frac{\partial x_{k+1}}{\partial x_k} \\
\frac{\partial \psi_{k+1}}{\partial x_k} = \frac{\partial \psi_{k+1}}{\partial x_{k+1}} \frac{\partial f_k}{\partial x_k}
\end{gather}
and the results for the controlled dynamical system $x_{k+1} = f\left(x_k,u_k\right)$ follow analogously by partially differentiating $\psi\left(x_{k+1}\right)$ with respect to the state and control inputs, respectively.
\end{proof}

\subsection{Continuous Time}
\label{Continuous Time}

In this section, we focus on two classes of Koopman representations, which we introduce and analyze for dynamical consistency in the following results. Both suppose that the Koopman operators are independent of the state and control inputs. Where they differ is in whether the state and control inputs are independently or jointly observed for the purposes of computing the control response.
\vspace{.5pc}
\begin{theorem} \label{thm:contsep}
If the Koopman representation of the form
\begin{gather} \label{eqn:koopcontsep}
\frac{\partial }{\partial t} \left(\psi_x\left(x\right)\right) = L_x \psi_x \left(x\right) + L_u \psi_u \left(u\right)
\end{gather}
is consistent with a controlled continuous-time dynamical system $\dot{x} = f\left( x, u \right)$ of the form
\begin{gather}
f\left(x,u\right) = f_x\left(x\right) + f_u \left(u\right) + f_{xu} \left(x,u\right) \\
f_u \left(0\right) = f_{xu} \left(x,0\right) = f_{xu} \left(0,u\right) = 0 \\
\psi_u \left(0\right) = 0
\end{gather}
then the following conditions on the Koopman and dynamical systems are true:
\begin{gather}
\frac{\partial \psi_x}{\partial x} f_x\left(x\right) = L_x \psi_x \left(x\right) \label{eqn:thmcontsep1}\\
\left. \frac{\partial \psi_x}{\partial x}\right|_{x=0} f_u\left(u\right) = L_u \psi_u \left(u\right) \label{eqn:thmcontsep2}\\
\left(\frac{\partial \psi_x}{\partial x} - \left. \frac{\partial \psi_x}{\partial x}\right|_{x=0}\right) f_u\left(u\right) \nonumber \\
+ \frac{\partial \psi_x}{\partial x} f_{xu} \left(x,u\right) = 0 \label{eqn:thmcontsep3}
\end{gather}
\end{theorem}
\vspace{.5pc}
\begin{proof}
Given the Koopman representation (\ref{eqn:koopcontsep}) of the continuous dynamical system with control defined in Theorem \ref{thm:contsep}, then dynamical consistency (viz. Theorem \ref{thm:dynconsist}) implies
\begin{align}
\frac{\partial }{\partial t} \left(\psi_x\left(x\right)\right) &= L_x \psi_x \left(x\right) + L_u \psi_u \left(u\right) \nonumber \\
&\hspace{-5ex}= \frac{\partial \psi_x}{\partial x} f\left(x,u\right) \nonumber \\
&\hspace{-5ex}= \frac{\partial \psi_x}{\partial x} f_x\left(x\right) + \frac{\partial \psi_x}{\partial x} f_u\left(u\right) + \frac{\partial \psi_x}{\partial x} f_{xu} \left(x,u\right) \label{eqn:kcs_expand}
\end{align}
Evalulating this at $u=0$ and $x=0$, respectively, yields the conditions
\begin{gather}
L_x \psi_x \left(x\right) = \frac{\partial \psi_x}{\partial x} f_x\left(x\right) \label{eqn:kcs_cond1} \\
L_u \psi_u \left(u\right) = \left. \frac{\partial \psi_x}{\partial x}\right|_{x=0} f_u\left(u\right) \label{eqn:kcs_cond2}
\end{gather}
Then, substituting (\ref{eqn:kcs_cond1}) and (\ref{eqn:kcs_cond2}) into (\ref{eqn:kcs_expand}) yields the last condition,
\begin{equation} \label{eqn:kcs_cond3}
\left(\frac{\partial \psi_x}{\partial x} - \left. \frac{\partial \psi_x}{\partial x}\right|_{x=0}\right) f_u\left(u\right) + \frac{\partial \psi_x}{\partial x} f_{xu} \left(x,u\right) = 0
\end{equation}
\end{proof}
\vspace{.5pc}
The following corollaries apply Theorem \ref{thm:contsep} to evaluate the generality of the continuous-time version of the representations proposed by Korda and Mezi\'{c} \cite{korda18jsr} and Arbabi et al. \cite{arbabi18jsr}.
\vspace{.5pc}
%
%
%
\begin{corollary} \label{cor:contsep1}
If the Koopman formulation in (\ref{eqn:koopcontsep}) is dynamically consistent and $\psi_x$ is state-inclusive, then $f\left( x, u \right) = 0$.
\end{corollary}
\vspace{.5pc}
\begin{proof}
Since $\psi_x$ is state-inclusive, then it includes $\mbox{Id}[x] = x$. Since the Jacobian of this subset is simply the identity matrix (and therefore constant), evaluating (\ref{eqn:kcs_cond3}) for this subset of observables yields the result $f_{xu}\left( x, u \right) = 0$.
\end{proof}
\vspace{.5pc}

By Corollary \ref{cor:contsep1}, if $\psi_x$ is state-inclusive and $f_{xu} \left(x,u\right)\neq 0$, then the formulation in (\ref{eqn:koopcontsep}) is not consistent with the controlled dynamical system. State-inclusivity limits the types of controlled dynamical systems for which separable Koopman representations of the form (\ref{eqn:koopcontsep}) are consistent.
\vspace{.5pc}


%

\begin{corollary} \label{cor:contsep2}
If the Koopman representation (\ref{eqn:koopcontsep}) is dynamically consistent, and if $f_{xu} \left(x,u\right) = 0$, then
\begin{gather}
\left( \left.\frac{\partial \psi_x}{\partial x}\right|_{x = x_1} - \left.\frac{\partial \psi_x}{\partial x}\right|_{x = x_2} \right) f_u\left(u\right) \nonumber \\
= 0 \ \forall \ x_1,x_2,u
\end{gather}
\end{corollary}
\vspace{.5pc}
\begin{proof}
This is a direct extension of condition (\ref{eqn:kcs_cond3}) in Theorem \ref{thm:contsep}:

\begin{gather}
\left(\left.\frac{\partial \psi_x}{\partial x}\right|_{x=x_1} - \left. \frac{\partial \psi_x}{\partial x}\right|_{x=0}\right) f_u\left(u\right) = 0 \nonumber \\
= \left(\left.\frac{\partial \psi_x}{\partial x}\right|_{x=x_2} - \left. \frac{\partial \psi_x}{\partial x}\right|_{x=0}\right) f_u\left(u\right) \ \forall \ x_1,x_2,u \\
\Rightarrow \left( \left.\frac{\partial \psi_x}{\partial x}\right|_{x = x_1} - \left.\frac{\partial \psi_x}{\partial x}\right|_{x = x_2} \right) f_u\left(u\right) \nonumber \\
= 0 \ \forall \ x_1,x_2,u
\end{gather}
\end{proof}
\vspace{.5pc}

This necessary consistency condition is not impossible to satisfy, but it is rather restrictive.  An easy way to satisfy it would be to restrict $\psi_x\left(x\right)$ to the identity observable, $\mbox{ID}[x]=x$, so that $\frac{\partial \psi_x}{\partial x}$ is constant, but this set of observables is not very expressive.
\vspace{.5pc}
\begin{corollary}
Consider the continuous-time version of the Koopman formulation in (\ref{eqn:kma_form}),
\begin{equation} \label{eqn:kma_cont}
\frac{\partial }{\partial t} \left(\psi_x \left(x\right) \right) = L \psi_x \left(x\right) + B u
\end{equation}
If $\psi_x$ is state-inclusive and (\ref{eqn:kma_cont}) is dynamically consistent, then the following conditions are true:
\begin{gather}
f_{xu} \left(x,u\right) = 0 \label{eqn:kma_cor1}\\
\left( \left.\frac{\partial \psi_x}{\partial x}\right|_{x = x_1} - \left.\frac{\partial \psi_x}{\partial x}\right|_{x = x_2} \right) f_u\left(u\right) \nonumber \\
= 0 \ \forall \ x_1,x_2,u \label{eqn:kma_cor2} \\
\left. \frac{\partial \psi_x}{\partial x} \right|_{x=0} \frac{\partial f_u}{\partial u} = B \ \forall \ x,u \label{eqn:kma_cor3} \\
\frac{\partial \psi_x}{\partial x} f_x\left(x\right) = L \psi_x \left(x\right) \label{eqn:kma_cor4}
\end{gather}

\end{corollary}
\vspace{.5pc}
\begin{proof}
Condition (\ref{eqn:kma_cor1}) follows from Corollary \ref{cor:contsep1}, just as condition (\ref{eqn:kma_cor2}) follows from Corollary \ref{cor:contsep2}. The proof of condition (\ref{eqn:kma_cor3}) follows directly from applying (\ref{eqn:thmcontsep2}) from Theorem \ref{thm:contsep} to the Koopman formulation in (\ref{eqn:kma_cont}):
\begin{gather}
L_u \psi_u \left(u\right) = Bu = \left. \frac{\partial \psi}{\partial x} \right|_{x=0} \cdot f_u\left(u\right) \\
\frac{\partial }{\partial u} \left(Bu\right) = \frac{\partial }{\partial u} \left( \left.\frac{\partial \psi}{\partial x}\right|_{x=0} f_u\left(u\right)\right) \\
B = \left.\frac{\partial \psi}{\partial x}\right|_{x=0} \frac{\partial f_u}{\partial u}
\end{gather}
Finally, the condition in (\ref{eqn:kma_cor4}) follows from Theorem \ref{thm:contsep} after setting $u=0$ in (\ref{eqn:kma_cont}).
\end{proof}
\vspace{.5pc}
It would be possible to satisfy these conditions if each observable $\psi\left(x\right)$ was linear in $x$ and $f_u\left(u\right)$ was linear in $u$; these are very restrictive conditions on the dynamical system and observables. In the following theorem, we consider an alternative Koopman representation where the control response is determined by jointly observing the state and control inputs.
\vspace{.5pc}

\begin{theorem} \label{thm:contjoint}
If the Koopman representation
\begin{equation} \label{eqn:contjoint}
\frac{\partial }{\partial t} \left(\psi_x\left(x\right)\right) = L_x \psi_x\left(x\right) + L_{xu} \psi_{xu} \left(x,u\right)
\end{equation}
is consistent with respect to the controlled continuous dynamical system $\dot{x} = f\left(x,u\right)$, such that
\begin{gather}
\psi_{xu} \left(x,0\right) = 0, \label{eqn:contjoint1} \\
f\left(x,u\right) = f_x \left(x\right) + f_{xu} \left(x,u\right), \label{eqn:contjoint2} \\
f_{xu}\left(x,0\right) = 0 \label{eqn:contjoint3}
\end{gather}
then the following conditions on the Koopman operators, observables, and dynamics are true:
\begin{gather}
\frac{\partial \psi_x}{\partial x} f_x\left(x\right) = L_x \psi_x\left(x\right) \\
\frac{\partial \psi_x}{\partial x} f_{xu}\left(x,u\right) = L_{xu} \psi_{xu} \left(x,u\right)
\end{gather}
\end{theorem}
\vspace{.5pc}
\begin{proof}
If the Koopman formulation (\ref{eqn:contjoint}) is consistent, then
\begin{equation}
L_x \psi_x\left(x\right) + L_{xu} \psi_{xu} \left(x,u\right) = \frac{\partial \psi_x}{\partial x} f\left(x,u\right)
\end{equation}
For $u=0$, 
\begin{gather}
L_x \psi_x\left(x\right) + L_{xu} \psi_{xu} \left(x,0\right) = \frac{\partial \psi_x}{\partial x} f\left(x,0\right) \\
L_{xu} \left( \psi_{xu} \left(x,0\right) - \psi_{xu} \left(x,0\right) \right) \nonumber \\
= \frac{\partial \psi_x}{\partial x} \left(f\left(x,u\right) - f\left(x,0\right)\right)
\end{gather}
Now, given the conditions of Theorem \ref{thm:contjoint}, $\psi_{xu} \left(x,0\right) = 0$, $f\left(x,u\right) = f_x \left(x\right) + f_{xu} \left(x,u\right)$, and $f_{xu}\left(x,0\right) = 0$, then the above simpifies to the theorem results.
\end{proof}
\vspace{.5pc}

This similar result to Theorem \ref{thm:contsep} shows a cleaner separation of terms than was possible for the separable Koopman representation in (\ref{eqn:koopcontsep}). It permits a state-inclusive set of observables without requiring simplifications to the dynamics in exchange for a larger set of control response observables, $\psi_{xu}(x,u)$.

The formulation of Kaiser et al. \cite{kaiser17jsr} given in (\ref{eqn:kaiser_form}) is also capable of being consistent as long as

\begin{gather}
\frac{\partial \psi}{\partial x} f\left(x,u\right) = \Lambda \psi \left(x,u\right)
\end{gather}

\noindent which seems straightforward and relatively unrestrictive.  As Kaiser et al. note, though, implementing this formulation requires controlling $\dot{u}$ rather than directly controlling $u$.

\subsection{Discrete Time \label{subsec:disc}}

In this section, we focus on two classes of discrete Koopman representations, which we introduce and analyze for dynamical consistency in the following results. Both discrete representations assume the Koopman operators are constant with respect to the state and control inputs. However, in one case the control response is dependent only on the control input, while in the other case the response is driven by joint observations of the state and control inputs.
\vspace{.5pc}
\begin{theorem} \label{thm:discsep}
If the Koopman representation
\begin{equation} \label{eqn:koopdiscsep}
\psi_x\left(x_{k+1}\right) = K_x \psi_x \left(x_k\right) + K_u \psi_u \left(u_k\right)
\end{equation}
is consistent with respect to the discrete-time dynamical system with control, $x_{k+1} = f\left(x_k,u_k\right)$, such that
\begin{gather}
f\left(x,u\right) = f_x\left(x\right) + f_u \left(u\right) + f_{xu} \left(x,u\right) \\
f_u \left(0\right) = f_{xu} \left(x,0\right) = f_{xu} \left(0,u\right) = 0 \\
\psi_u \left(0\right) = 0
\end{gather}
then the following conditions on the Koopman and dynamical systems are true:
\begin{gather}
\left.\frac{\partial \psi_{x,k+1}}{\partial x_{k+1}}\right|_{u_k = 0} \frac{\partial f_{x,k}}{\partial x_k} = K_x \frac{\partial \psi_{x,k}}{\partial x_k} \label{eqn:discsepcond1}\\
\left. \frac{\partial \psi_{x,k+1}}{\partial x_{k+1}}\right|_{x_k=0} \frac{\partial f_{u,k}}{\partial u_k} = K_u \frac{\partial \psi_{u,k}}{\partial u_k} \label{eqn:discsepcond2}\\
\left(\frac{\partial \psi_{x,k+1}}{\partial x_{k+1}} - \left. \frac{\partial \psi_{x,k+1}}{\partial x_{k+1}}\right|_{x_k=0}\right) \frac{\partial f_{u,k}}{\partial u_k} \nonumber \\
+ \frac{\partial \psi_{u,k}}{\partial u_k} \frac{\partial f_{xu,k}}{\partial u_k} = 0 \label{eqn:discsepcond3}\\
\left(\frac{\partial \psi_{x,k+1}}{\partial x_{k+1}} - \left. \frac{\partial \psi_{x,k+1}}{\partial x_{k+1}}\right|_{u_k=0}\right) \frac{\partial f_{x,k}}{\partial x_k} \nonumber \\
+ \frac{\partial \psi_{x,k}}{\partial x_k} \frac{\partial f_{xu,k}}{\partial x_k} = 0 \label{eqn:discsepcond4}
\end{gather}
\end{theorem}
\vspace{.5pc}
\begin{proof}
If the Koopman representation (\ref{eqn:koopdiscsep}) is consistent, then it can be shown that
\begin{gather}
\frac{\partial \psi_{x,k+1}}{\partial x_k} = \frac{\partial \psi_{x,k+1}}{\partial x_{k+1}} \frac{\partial f_k}{\partial x_k} = K_x \frac{\partial \psi_{x,k}}{\partial x_k} \\
\frac{\partial \psi_{x,k+1}}{\partial x_{k+1}} \frac{\partial f_{x,k}}{\partial x_k} + \frac{\partial \psi_{x,k+1}}{\partial x_{k+1}} \frac{\partial f_{xu,k}}{\partial x_k} = K_x \frac{\partial \psi_{x,k}}{\partial x_k} \\
\frac{\partial \psi_{x,k+1}}{\partial u_k} = \frac{\partial \psi_{x,k+1}}{\partial x_{k+1}} \frac{\partial f_k}{\partial u_k} = K_u \frac{\partial \psi_{u,k}}{\partial u_k} \\
\frac{\partial \psi_{x,k+1}}{\partial x_{k+1}} \frac{\partial f_{u,k}}{\partial u_k} + \frac{\partial \psi_{x,k+1}}{\partial x_{k+1}} \frac{\partial f_{xu,k}}{\partial u_k} = K_x \frac{\partial \psi_{u,k}}{\partial u_k}
\end{gather}
Evaluating at $x_k = 0$ yields
\begin{gather}
\left. \frac{\partial \psi_{x,k+1}}{\partial x_{k+1}}\right|_{x_k=0} \frac{\partial f_{u,k}}{\partial u_k} = K_u \frac{\partial \psi_{u,k}}{\partial u_k} \\
\left(\frac{\partial \psi_{x,k+1}}{\partial x_{k+1}} - \left. \frac{\partial \psi_{x,k+1}}{\partial x_{k+1}}\right|_{x_k=0}\right) \frac{\partial f_{u,k}}{\partial u_k} \nonumber \\
+ \frac{\partial \psi_{u,k}}{\partial u_k} \frac{\partial f_{xu,k}}{\partial u_k} = 0
\end{gather}
and evaluating at $u_k = 0$ yields
\begin{gather}
\left. \frac{\partial \psi_{x,k+1}}{\partial x_{k+1}}\right|_{x_k=0} \frac{\partial f_{u,k}}{\partial u_k} = K_u \frac{\partial \psi_{u,k}}{\partial u_k} \\
\left(\frac{\partial \psi_{x,k+1}}{\partial x_{k+1}} - \left. \frac{\partial \psi_{x,k+1}}{\partial x_{k+1}}\right|_{u_k=0}\right) \frac{\partial f_{x,k}}{\partial x_k} \nonumber \\
+ \frac{\partial \psi_{x,k}}{\partial x_k} \frac{\partial f_{xu,k}}{\partial x_k} = 0
\end{gather}
Since $f_{xu} \left(x,0\right) = f_{xu} \left(0,u\right) = 0$,

\begin{gather}
\left. \frac{\partial f_{xu}}{\partial x} \right|_{u=0} = \left. \frac{\partial f_{xu}}{\partial u} \right|_{x=0} = 0
\end{gather}

\end{proof}
\vspace{.5pc}

\begin{corollary} \label{cor:discsep1}
If the Koopman representation (\ref{eqn:koopdiscsep}) is consistent and the state observables $\psi_x$ are state-inclusive, then $f_{xu}\left(x_k, u_k\right) = 0$.
\end{corollary}
\vspace{.5pc}
\begin{proof}
Since (\ref{eqn:koopdiscsep}) is consistent, then (\ref{eqn:discsepcond3}) and (\ref{eqn:discsepcond4}) are true, and if $\psi_x$ is state inclusive, then evaluating (\ref{eqn:discsepcond3}) and (\ref{eqn:discsepcond3}) for the subset of observables $\mbox{Id}[x] = x$ yields
\begin{gather}
\frac{\partial f_{xu,k}}{\partial u_k} = 0\\
\frac{\partial f_{xu,k}}{\partial x_k} = 0
\end{gather}

Therefore $f_{xu} \left(x,u\right)$ is constant, and since $f_{xu}\left(x,0\right) = f_{xu} \left(0,u\right) = 0$, $f_{xu} \left(x,u\right) = 0$.
\end{proof}
\vspace{.5pc}

\begin{corollary} \label{cor:discsep2}
If the Koopman representation (\ref{eqn:koopdiscsep}) is consistent, and if $f_{xu} \left(x,u\right) = 0$, then the system in Theorem \ref{thm:discsep} can only be consistent if
\begin{gather}
\left(\!\! \left.\frac{\partial \psi_{x,k+1}}{\partial x_{k+1}}\right|_{\overset{x_k=x_1 }{ u_k=u_1}} \!\! - \!\! \left. \frac{\partial \psi_{x,k+1}}{\partial x_{k+1}}\right|_{\overset{x_k=x_2 }{u_k=u_1}}\right)\!\! \left.\frac{\partial f_{u,k}}{\partial u_k}\right|_{u_k = u_1} \!\! = 0 \\
\left(\!\! \left.\frac{\partial \psi_{x,k+1}}{\partial x_{k+1}}\right|_{\overset{ x_k=x_1}{ u_k=u_1}} \!\! - \!\! \left. \frac{\partial \psi_{x,k+1}}{\partial x_{k+1}}\right|_{ \overset{x_k=x_1}{ u_k=u_2} }\right) \!\! \left. \frac{\partial f_{x,k}}{\partial x_k}\right|_{x_k=x_1} \!\! = 0
\end{gather}
\noindent for all $x_1,x_2,u_1,u_2$.
\end{corollary}
\vspace{.5pc}
\begin{proof}
This follows directly from Theorem \ref{thm:discsep}. If $f_{xu} \left(x,u\right) = 0$, then

\begin{gather}
\left(\left.\frac{\partial \psi_{x,k+1}}{\partial x_{k+1}}\right|_{x_k=x_1} - \left. \frac{\partial \psi_{x,k+1}}{\partial x_{k+1}}\right|_{x_k=0}\right) \frac{\partial f_{u,k}}{\partial u_k} = 0 \nonumber \\
= \left(\left.\frac{\partial \psi_{x,k+1}}{\partial x_{k+1}}\right|_{x_k=x_2} - \left. \frac{\partial \psi_{x,k+1}}{\partial x_{k+1}}\right|_{x_k=0}\right) \frac{\partial f_{u,k}}{\partial u_k} \\
\Rightarrow \left(\left.\frac{\partial \psi_{x,k+1}}{\partial x_{k+1}}\right|_{x_k=x_1} - \left. \frac{\partial \psi_{x,k+1}}{\partial x_{k+1}}\right|_{x_k=x_2}\right) \frac{\partial f_{u,k}}{\partial u_k} \nonumber \\
= 0 \ \forall \ x_1,x_2, u\\
\left(\left.\frac{\partial \psi_{x,k+1}}{\partial x_{k+1}}\right|_{u_k=u_1} - \left. \frac{\partial \psi_{x,k+1}}{\partial x_{k+1}}\right|_{u_k=0}\right) \frac{\partial f_{x,k}}{\partial x_k} = 0  \nonumber \\
= \left(\left.\frac{\partial \psi_{x,k+1}}{\partial x_{k+1}}\right|_{u_k=u_2} - \left. \frac{\partial \psi_{x,k+1}}{\partial x_{k+1}}\right|_{u_k=0}\right) \frac{\partial f_{x,k}}{\partial x_k} \\
\Rightarrow \left(\left.\frac{\partial \psi_{x,k+1}}{\partial x_{k+1}}\right|_{u_k=u_1} - \left. \frac{\partial \psi_{x,k+1}}{\partial x_{k+1}}\right|_{u_k=u_2}\right) \frac{\partial f_{x,k}}{\partial x_k} \nonumber \\
= 0 \ \forall \ u_1,u_2,x
\end{gather}

\end{proof}
\vspace{.5pc}
\begin{corollary}
If the Koopman formulation $\psi \left(x_{k+1}\right) = K \psi \left(x_k\right) + B u_k$ is dynamically consistent and $\psi \left(x \right)$ is state-inclusive, then
\begin{gather}
f_{xu} \left(x,u\right) = 0 \\
\frac{\partial \psi_{x,k+1}}{\partial x_{k+1}} \frac{\partial f_{u,k}}{\partial u_k} =  B \ \forall \ x,u
\end{gather}
\end{corollary}
\vspace{.5pc}
\begin{proof}
The first condition is implied by Corollary \ref{cor:discsep1}.  Given that, consistency implies that
\begin{gather}
\frac{\partial }{\partial u_k} \left(B u_k\right) = \frac{\partial \psi_{x,k+1}}{\partial x_{k+1}} \frac{\partial f_{u,k}}{\partial u_k} \\
B = \frac{\partial \psi_{x,k+1}}{\partial x_{k+1}} \frac{\partial f_{u,k}}{\partial u_k}
\end{gather}
\end{proof}
\vspace{.5pc}

This is not impossible, but it does seem rather restrictive; the easiest way to satisfy it would be to have $\psi_x\left(x\right)$ as a linear function of $x$, so that $\frac{\partial \psi_x}{\partial x}$ is constant, and to have $f_u\left(u\right)$ to be a linear function of $u$.
\vspace{.5pc}
\begin{theorem}
If the Koopman representation
\begin{gather}
\psi_x\left(x_{k+1}\right) = K_x \psi_x\left(x_k\right) + K_{xu} \psi_{xu} \left(x_k,u_k\right) \label{eqn:discjoint}
\end{gather}
\noindent is consistent with respect to $x_{k+1}= f\left(x_k,u_k\right)$, then
\begin{gather}
\left[\frac{\partial \psi_{x,k+1}}{\partial x_{k+1}} \frac{\partial f_k}{\partial x_k}\right]_{u_k=0} = K_x \frac{\psi_{x,k}}{\partial x_k} \nonumber \\
+ K_{xu} \left.\frac{\partial \psi_{xu,k}}{\partial x_k}\right|_{u_k=0} \\
\frac{\partial \psi_{x,k+1}}{\partial x_{k+1}} \frac{\partial f_k}{\partial u_k} = K_{xu} \frac{\partial \psi_{xu,k}}{\partial u_k}
\end{gather}
\end{theorem}
\vspace{.5pc}
\begin{proof}
For $u_k=0$, consistency implies that
\begin{gather}
\frac{\partial \psi_{x,k+1}}{\partial x_k} = \frac{\partial \psi_{x,k+1}}{\partial x_{k+1}} \frac{\partial f_k}{\partial x_k} \\
K_x \frac{\psi_{x,k}}{\partial x_k} + K_{xu} \left.\frac{\partial \psi_{xu,k}}{\partial x_k}\right|_{u_k=0} \nonumber \\
= \left[\frac{\partial \psi_{x,k+1}}{\partial x_{k+1}} \frac{\partial f_k}{\partial x_k}\right]_{u_k=0}
\end{gather}
To be consistent when control is applied, i.e., $x_{k+1} = f\left(x_k,u_k\right)$, then
\begin{gather}
\frac{\partial \psi_{x,k+1}}{\partial u_k} = \frac{\partial \psi_{x,k+1}}{\partial x_{k+1}} \frac{\partial f_k}{\partial u_k} \\
K_{xu} \frac{\partial \psi_{xu,k}}{\partial u_k} = \frac{\partial \psi_{x,k+1}}{\partial x_{k+1}} \frac{\partial f_k}{\partial u_k}
\end{gather}
\end{proof}
\vspace{.5pc}
\begin{corollary}
If in addition it is assumed that $\psi_{xu} \left(x,0\right) = 0$, then the above conditions simplify down to
\begin{gather}
\left[\frac{\partial \psi_{x,k+1}}{\partial x_{k+1}} \frac{\partial f_k}{\partial x_k}\right]_{u_k=0} = K_x \frac{\psi_{x,k}}{\partial x_k} \\
\frac{\partial \psi_{x,k+1}}{\partial x_{k+1}} \frac{\partial f_k}{\partial u_k} = K_{xu} \frac{\partial \psi_{xu,k}}{\partial u_k}
\end{gather}
\end{corollary}
\vspace{.5pc}
\begin{corollary}
If it is further assumed that $\psi_{xu} \left(x,0\right) = 0$, and $f\left(x,u\right) = f_x \left(x\right) + f_{xu} \left(x,u\right)$ such that $f_{xu}\left(x,0\right) = 0$, then the above conditions further simplify to
\begin{gather}
\left[\frac{\partial \psi_{x,k+1}}{\partial x_{k+1}} \frac{\partial f_{x,k}}{\partial x_k}\right]_{u_k=0} = K_x \frac{\psi_{x,k}}{\partial x_k} 
\label{eqn:xu_consist1}\\
\frac{\partial \psi_{x,k+1}}{\partial x_{k+1}} \frac{\partial f_{xu,k}}{\partial u_k} = K_{xu} \frac{\partial \psi_{xu,k}}{\partial u_k}
\label{eqn:xu_consist2}
\end{gather}
\end{corollary}
\vspace{.5pc}
The representation of Williams et al. \cite{williams16cp} is directly related to (\ref{eqn:discjoint}). Recall that
\begin{gather}
\psi\left(x_{k+1}\right) = K\left(u_k\right) \psi \left(x_k\right) \\
K\left(u\right) = \sum_i \psi^u_i \left(u\right) K_i
\end{gather}
To turn this into a form amenable to linear control techniques, consider
\begin{gather}
\psi_x \left(x_{k+1}\right) = K\left(0\right) \psi_x \left(x_k\right) + I \psi_{xu} \left(x_k,u_k\right) \\
\psi_{xu} \left(x,u\right) = \left(K\left(u\right) - K\left(0\right)\right) \psi_x\left(x\right) 
\label{williams constr}
\end{gather}
The consistency conditions in (\ref{eqn:xu_consist1}) and (\ref{eqn:xu_consist2}) then become

\begin{gather}
\left[\frac{\partial \psi_{x,k+1}}{\partial x_{k+1}} \frac{\partial f_{x,k}}{\partial x_k}\right]_{u_k=0} = K\left(0\right) \frac{\psi_{k}}{\partial x_k} \\
\frac{\partial \psi_{x,k+1}}{\partial x_{k+1}} \frac{\partial f_{xu,k}}{\partial u_k} = \left[\sum_i \frac{\partial \psi^u_i}{\partial u} K_i \right]_{u=u_k} \psi \left(x_k\right)
\end{gather}

\section{Discussion}
\label{discussion}
Thus far, there has not been previous work on whether particular Koopman operator control formulations are actually capable, in principle, of representing systems of interest.  It is worth knowing when it is possible to use certain formulations. Using (\ref{eqn:koopcontsep}) and simplifications thereof is attractive in that it makes it possible to separate state and control variables.  However, the results in the previous section show that this formulation comes with some limitations -- especially if $\psi_x\left(x\right)$ is state-inclusive.  If $\psi_x\left(x\right)$ is state-inclusive, then it is necessary to find observables that essentially solve a system of equations of the form
\begin{gather}
G\left(x\right) f\left(u\right) = 0 \ \forall \ x,u
\end{gather}
where $G\left(x\right)$ is a matrix and the variable of interest, and $f\left(u\right)$ is a known vector function.  This system can be solved if any non-zero $f\left(u\right)$ remains in the nullspace of $G\left(x\right)$.  It is a strong constraint on $G\left(x\right)$ (and $f\left(u\right)$, though.  In the context of Koopman operator control, this system functions like an additional linearity enforced in the system -- one observed along trajectories -- on top of the linear dynamics of the observables.  It is not clear how this constraint interacts with the need for the observables to be invariant.

There are Koopman operator control formulations that are less constrained.  The above results show that the formulation of Williams et al. \cite{williams16cp} is essentially a special case of (\ref{eqn:contjoint}), both of which are more general.  In particular, they can model a nonlinear system with state-dependent control (i.e. $f_{xu}\left(x,u\right) \neq 0$) while still having state-inclusive observables.

The formulation in (\ref{eqn:contjoint}) and its discrete-time counterpart also allow for a cleaner distinction between the controlled and uncontrolled dynamics.  In (\ref{eqn:koopcontsep}), for example, solving for $\psi_x$ by considering the uncontrolled system would only lead to satisfying (\ref{eqn:thmcontsep1}).  However, to be consistent for the control system as well, it is also necessary to satisfy (\ref{eqn:thmcontsep3}) for all $u$, and if $\psi_x$ is already determined, all of the quantities in (\ref{eqn:thmcontsep3}) equation are fixed.  For (\ref{eqn:contjoint}), however, by imposing $\psi_{xu} \left(x,0\right) = 0$, which seems reasonable, and then solving for the uncontrolled system satisfies (\ref{eqn:xu_consist1}), one can then calculate the Koopman representation for the controlled system by taking $\psi_x$ as fixed and solving (\ref{eqn:xu_consist2}) for $\psi_{xu}$.

In practice, data-driven finite numerical approximations to the Koopman operator always have some error associated with them -- they are not perfectly invariant, and they do not exactly follow the dynamics of the underlying nonlinear system.  Using (\ref{eqn:koopcontsep}) may be more computationally efficient than (\ref{eqn:contjoint}) when applying linear control methods, and any differences in fidelity to the original dynamical system may be worthwhile to accelerate computation.  However, those seeking to apply control techniques to Koopman operators should be aware of the existence and nature of those tradeoffs.

\section{Conclusion}

Koopman operators are used to simplify the analysis of nonlinear dynamical systems by embedding their flows in a (possibly linear) space of observables. Koopman representations have not been applied widely to dynamical systems with control. In this paper, several Koopman formulations were analyzed and compared in terms of the necessary conditions imposed on the dynamical system and Koopman parameters. A consistency condition was introduced for determining necessary conditions for a Koopman representation to admit a given dynamical system, or otherwise introduce model error. By both independently and jointly observing the state and control inputs, the resulting formulation was shown to admit a relatively wide class of dynamical systems.

\section{ACKNOWLEDGMENTS}

The research described in this paper is part of the Agile Deep Science Initiative at Pacific Northwest National Laboratory. It was conducted under the Laboratory Directed Research and Development Program at PNNL, a multiprogram national laboratory operated by Battelle for the U.S. Department of Energy.

%
%


\bibliography{bib}

\begin{thebibliography}{10}

\bibitem{budisic12jsr}
M.~Budi{\v{s}}i{\'c}, R.~Mohr, and I.~Mezi{\'c}, ``Applied koopmanism,'' {\em
  Chaos: An Interdisciplinary Journal of Nonlinear Science}, vol.~22, no.~4,
  p.~047510, 2012.

\bibitem{page18jsr}
J.~Page and R.~R. Kerswell, ``Koopman analysis of burgers equation,'' {\em
  Physical Review Fluids}, vol.~3, no.~7, p.~071901, 2018.

\bibitem{tu14jsr}
J.~H. Tu, C.~W. Rowley, D.~M. Luchtenburg, S.~L. Brunton, and J.~N. Kutz, ``On
  dynamic mode decomposition: Theory and applications,'' {\em Journal of
  Computational Dynamics}, vol.~1, no.~2, pp.~391--421, 2014.

\bibitem{williams15jsr}
M.~O. Williams, I.~G. Kevrekidis, and C.~W. Rowley, ``A data--driven
  approximation of the koopman operator: Extending dynamic mode
  decomposition,'' {\em Journal of Nonlinear Science}, vol.~25, no.~6,
  pp.~1307--1346, 2015.

\bibitem{kutz16jsr}
J.~N. Kutz, J.~L. Proctor, and S.~L. Brunton, ``Koopman theory for partial
  differential equations,'' {\em arXiv preprint arXiv:1607.07076}, 2016.

\bibitem{proctor18jsr}
J.~L. Proctor, S.~L. Brunton, and J.~N. Kutz, ``Generalizing koopman theory to
  allow for inputs and control,'' {\em SIAM Journal on Applied Dynamical
  Systems}, vol.~17, no.~1, pp.~909--930, 2018.

\bibitem{williams16cp}
M.~O. Williams, M.~S. Hemati, S.~T. Dawson, I.~G. Kevrekidis, and C.~W. Rowley,
  ``Extending data-driven koopman analysis to actuated systems,'' in {\em IFAC
  Symposium on Nonlinear Control Systems (NOLCOS)}, 2016.

\bibitem{korda18jsr}
M.~Korda and I.~Mezi{\'c}, ``Linear predictors for nonlinear dynamical systems:
  Koopman operator meets model predictive control,'' {\em Automatica}, vol.~93,
  pp.~149--160, 2018.

\bibitem{williams14jsr}
M.~O. Williams, C.~W. Rowley, and I.~G. Kevrekidis, ``A kernel-based approach
  to data-driven koopman spectral analysis,'' {\em arXiv preprint
  arXiv:1411.2260}, 2014.

\bibitem{yeung17jsr}
E.~Yeung, S.~Kundu, and N.~Hodas, ``Learning deep neural network
  representations for koopman operators of nonlinear dynamical systems,'' {\em
  arXiv preprint arXiv:1708.06850}, 2017.

\bibitem{li17jsr}
Q.~Li, F.~Dietrich, E.~M. Bollt, and I.~G. Kevrekidis, ``Extended dynamic mode
  decomposition with dictionary learning: A data-driven adaptive spectral
  decomposition of the koopman operator,'' {\em Chaos: An Interdisciplinary
  Journal of Nonlinear Science}, vol.~27, no.~10, p.~103111, 2017.

\bibitem{brunton16jsr}
S.~L. Brunton, B.~W. Brunton, J.~L. Proctor, and J.~N. Kutz, ``Koopman
  invariant subspaces and finite linear representations of nonlinear dynamical
  systems for control,'' {\em PloS one}, vol.~11, no.~2, p.~e0150171, 2016.

\bibitem{proctor16jsr}
J.~L. Proctor, S.~L. Brunton, and J.~N. Kutz, ``Dynamic mode decomposition with
  control,'' {\em SIAM Journal on Applied Dynamical Systems}, vol.~15, no.~1,
  pp.~142--161, 2016.

\bibitem{kaiser17jsr}
E.~Kaiser, J.~N. Kutz, and S.~L. Brunton, ``Data-driven discovery of koopman
  eigenfunctions for control,'' {\em arXiv preprint arXiv:1707.01146}, 2017.

\bibitem{arbabi18jsr}
H.~Arbabi, M.~Korda, and I.~Mezic, ``A data-driven koopman model predictive
  control framework for nonlinear flows,'' {\em arXiv preprint
  arXiv:1804.05291}, 2018.

\end{thebibliography}
\bibliographystyle{ieeetr}

\end{document}